\documentclass[a4paper]{article}
\usepackage{geometry}
\usepackage{amssymb}
\usepackage{amsmath,cases}
\usepackage{amsthm}
\usepackage{enumerate}
\usepackage{xcolor}
\usepackage{graphicx}
\usepackage{textcomp}
\usepackage{mathrsfs}
\usepackage{bbm}
\usepackage{bm}
\usepackage{caption}
\usepackage{chemarrow}
\usepackage{stmaryrd}
\usepackage{url}
\usepackage{psfrag}
\usepackage[active]{srcltx}   
\usepackage[all,cmtip]{xy}
\usepackage{pdfsync}																				
\usepackage{yhmath}
\usepackage[bookmarksnumbered=true, pdfauthor={Jérémie Bettinelli}]{hyperref}									

\newcommand{\choo}[2]{\begin{pmatrix}#1 \\ #2 \end{pmatrix}}

\newcommand{\de}{\mathrel{\mathop:}\hspace*{-.6pt}=}

\newcommand{\ooo}{\bullet}

\newcommand{\F}{\mathcal{F}}

\newcommand{\Q}{\mathcal{Q}}

\newcommand{\Qnp}{\mathcal{Q}_{n,p}}

\newcommand{\f}{\mathfrak{f}}

\newcommand{\lab}{\mathfrak{l}}

\newcommand{\m}{\mathfrak{m}}

\newcommand{\q}{\mathfrak{q}}

\newcommand{\tr}{\mathfrak{t}}

\newcommand{\BDG}{Bouttier--Di~Francesco--Guitter }

\newcommand{\rev}{\operatorname{rev}}
\newcommand{\up}{\text{\raisebox{.6pt}{\scriptsize $\uparrow$}}}
\newcommand{\dn}{\text{\raisebox{.6pt}{\scriptsize $\downarrow$}}}
\newcommand{\pc}{\mathbbm c}

\newcommand{\ph}{\mathbbm h}
\newcommand{\pl}{\mathbbm l}
\newcommand{\pp}{\mathbbm p}
\newcommand{\pr}{\mathbbm r}

\newcommand{\sand}{\qquad\text{ and }\qquad}

\newcommand{\lp}{\left(}
\newcommand{\rp}{\right)}

\definecolor{gris}{gray}{0.7}
\definecolor{grisf}{gray}{0.4}
\definecolor{vert}{rgb}{0,.5547,0}

\theoremstyle{plain}
\newtheorem{thm}{Theorem}

\newtheorem{defi}{Definition}

\theoremstyle{definition}
\newtheorem*{rem}{Remark}

\theoremstyle{remark}

\newenvironment{pre}[1][\proofname]{%
  \proof[\bfseries #1]%
}{\endproof}

\title{Increasing Forests and Quadrangulations via a Bijective Approach\thanks{This work is partially supported by ANR-08-BLAN-0190}}
\author{J\'er\'emie \textsc{Bettinelli}\\
	\small Institut \'Elie Cartan de Lorraine, B.P.\ 239\\
	\small F-54506 Vand\oe uvre-l\`es-Nancy Cedex\\
	\small\url{jeremie.bettinelli@normalesup.org}\\
	\small\url{http://www.normalesup.org/~bettinel}}

\begin{document}
\maketitle

\begin{abstract}

In this work, we expose four bijections each allowing to increase (or decrease) one parameter in either uniform random forests with a fixed number of edges and trees, or quadrangulations with a boundary having a fixed number of faces and a fixed boundary length. In particular, this gives a way to sample a uniform quadrangulation with $n+1$ faces from a uniform quadrangulation with~$n$ faces or a uniform forest with $n+1$ edges and~$p$ trees from a uniform forest with~$n$ edges and~$p$ trees.

\bigskip



\end{abstract}

\tableofcontents

\section{Introduction}

Maps are known to have a lot of applications in different fields of mathematics, computer science and physics. These applications strongly rely on their combinatorial structures, which have been widely investigated during the last few decades. A particularly interesting class of maps is the class of quadrangulations with or without a boundary, which are for example natural candidates for discretizing surfaces.

Our interest in this paper, inspired from R\'emy's algorithm~\cite{remy85} on growing trees, is in finding a natural way to grow a planar quadrangulation that preserves the uniform measure. For instance, we will present a bijection between quadrangulations with a boundary having~$n$ faces and~$2p$ half-edges on the boundary carrying some distinguished elements and quadrangulations with a boundary having $n+1$ faces and~$2p$ half-edges on the boundary also carrying distinguished elements. Our bijection is designed in such a way that the number of possibilities for distinguishing the necessary elements only depends on the size and boundary length of the quadrangulations. Forgetting these distinguished elements, we obtain a way to sample a uniform quadrangulation with a boundary having a prescribed number of faces and boundary length from a uniform quadrangulation with a boundary having one less face. In other words, there exist some integer constants $c_{n,p}$ and $c'_{n,p}$ such that our construction provides a $c_{n,p}$-to-$c'_{n,p}$ mapping between the set of quadrangulations with a boundary having~$n$ faces and~$2p$ half-edges on the boundary and the set of quadrangulations with a boundary having $n+1$ faces and~$2p$ half-edges on the boundary.

In addition to the probabilistic point of view, such bijections also present a combinatorial interest as they provide an interpretation to some combinatorial identities. The bijections we present in this work interpret already known identities so that they actually provide alternate proofs for these identities. In the future, we hope that other similar bijections will allow to solve some open enumeration problems.

Our method of ``cut and glue'' bijections, which could informally be pictured as unbuttoning a shirt and buttoning it back incorrectly by putting each button into the hole that immediately follows the correct one, possesses a certain robustness and can be declined in many ways. We present here in details four such bijections by focusing on forests and quadrangulations with a boundary. In an upcoming work, we plan to present more bijections relying on the same idea. In particular, one of these bijections will allow to recover Tutte's formula~\cite{tutte62cs} counting the number of planar maps with~$n$ faces having prescribed degrees~$a_1$, \ldots, $a_n$ where at most two~$a_i$'s are odd numbers.

It has also been pointed to us that our method somehow recalls a work by Cori~\cite{cori75code} where he used a so-called transfer bijection roughly consisting in transferring one degree from a face to a neighboring face. Using a properly defined chain, this allows to transfer one degree from a face to any other face, step by step. We do not believe that his results are related to the present work but we think they are worth mentioning at this point.

\section{Setting and presentation of the results}

Recall that a planar map is an embedding of a finite connected graph (possibly with loops and multiple edges) into the two-dimensional sphere, considered up to direct homeomorphisms of the sphere. The faces of the map are the connected components of the complement of edges. We will call \emph{half-edge} an edge carrying one of its two possible orientations. We say that a half-edge~$h$ is \emph{incident} to a face~$f$ (or that~$f$ is incident to~$h$) if~$h$ belongs to the boundary of~$f$ and is oriented in such a way that~$f$ lies to its left. A \emph{corner} is an angular sector of a face delimited by two consecutive half-edges incident to the face. We will implicitly consider our maps to be rooted, which means that one corner of one face is distinguished. This distinguished corner will be called the \emph{root corner}.

A \emph{quadrangulation with a boundary} is a particular instance of planar map whose faces are all incident to exactly~$4$ half-edges, with the exception of the face containing the root corner, which may be of arbitrary even degree. The latter face will be called the \textit{external face}, whereas the other ones will be called \textit{internal faces}. To match this terminology, we will as often as possible draw the external face as the infinite component of the plane on our figures. As a result, note that, since the cyclic ordering of edges around a vertex is prescribed, it is unambiguous to speak of clockwise and counterclockwise order. The half-edges incident to the external face will constitute the \textit{boundary} of the map. Beware that we do not require here the boundary to be a simple curve.

\begin{defi}
For $n\ge 0$ and $p \ge 1$, we will denote by $\Qnp$ the set of all quadrangulations with a boundary having~$n$ internal faces and~$2p$ half-edges on the boundary.
\end{defi}

\begin{figure}[ht]
	\centering\includegraphics[height=38mm]{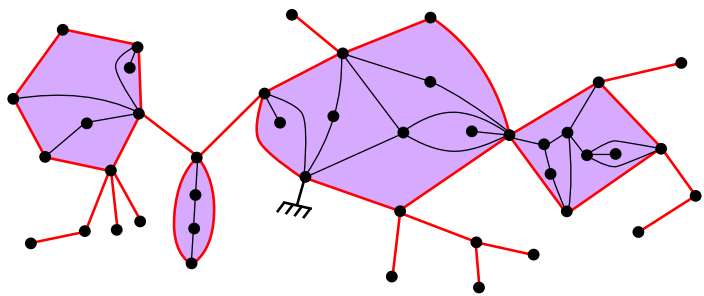}
	\caption{Example of quadrangulation with a boundary having~$19$ internal faces and~$46$ half-edges on the boundary. The root corner is represented by the rake-like symbol in the middle part of the picture.}
	\label{quad_bord}
\end{figure}

Let us make a few observations:
\begin{itemize}
	\item As we do not require the boundary of a quadrangulation with a boundary to be a simple curve, the set $\Qnp$ is never empty.
	\item There is a trivial bijection between $\Q_{n+1,1}$ and $\Q_{n,2}$: in the direct sense, just remove the one edge of the boundary that is directly to the right of the root. Conversely, double the edge directly to the left of the root in such a way that the degree-$2$ face so created contains the root. 
	\item The two previous sets correspond to the set of quadrangulations (without boundary) with $n+1$ faces, so our approach actually contains this case.
	\item The set $\Q_{0,p}$ is the set of rooted plane trees having~$p$ edges.
	\item All the cycles in a quadrangulation with a boundary have an even length. Indeed, if there were cycles with an odd length, let us consider one of these cycles. By Jordan curve theorem, it splits the map into two maps. Let us consider any of these maps. All its faces except one have an even degree. This is in contradiction with the fact that the number of half-edges is necessarily even.
	\item Using Euler's characteristic formula, it is easy to see that an element of~$\Qnp$ has exactly
		\begin{itemize}
			\item $n$ internal faces,
			\item $2p$ half-edges on the boundary,
			\item $n+p+1$ vertices,
			\item $2n+p$ edges.
		\end{itemize}
\end{itemize}

The cardinality of the set $\Qnp$ may be computed by several methods. We may use Tutte's formula~\cite{tutte62cs}, a generating series method~\cite{bouttier09dsqb} or a bijective approach using the \BDG bijection~\cite{bouttier04pml} (see for example~\cite{bettinelli11slr}):
\begin{equation}\label{qnp}
|\Qnp| =\frac{3^n\, (2p)!\, (2n+p-1)!}{p!\, (p-1)! \, n!\, (n+p+1)!}.
\end{equation}
In this work, we present two bijections accounting for the following combinatorial identities:
\begin{align}
(2p+1)(2p+2)(2n+p)\,|\Qnp| &= p\,(p+1)(n+p+2)\,|\Q_{n,p+1}|,\label{pp}\\
3(2n+p)(2n+p+1)\,|\Qnp|&=(n+1)(n+p+2)\,|\Q_{n+1,p}|.\label{nn}
\end{align}
In particular, using the initial condition $\Q_{0,1}=1$, our approach also provides a new proof for~\eqref{qnp}.

\bigskip

The strategy we use is roughly the following. Thanks to some distinguished elements of a map (faces, edges, vertices, corners), we construct a path in the map. Then, we ``cut'' along this path and ``glue'' back after shifting a little bit. This operation mildly modifies the map along the path and changes its structure around the extremities of the path, creating new distinguished elements. In order for this operation to work, the path we construct has to be totally recoverable from the new distinguished elements. In general, the notion of left-most geodesics or right-most geodesics will allow this.

As a warm-up, we present in Section~\ref{secf} two bijections between forests. Although not directly related to our bijections on quadrangulations, they use the same kind of ideas and are easier to handle. In addition, some definitions and notation that will be used throughout this paper are also provided at this point. We then present in Section~\ref{secpp} our bijection interpreting identity~\eqref{pp} and Section~\ref{secnn} is devoted to our last bijection, corresponding to~\eqref{nn}.

Let us also mention that our bijections on maps are not coming from bijections on forests that are transferred through the \BDG bijection~\cite{bouttier04pml}; our bijections work directly on the maps. See Section~\ref{secbdg} for more details regarding this matter. For more simplicity, we will from now on use the term quadrangulation to mean rooted quadrangulation with a boundary.

\section{Increasing forests}\label{secf}

This section presents two bijections allowing to grow uniform forests with a fixed number of trees and edges. These bijections are simpler than our main bijections for quadrangulations and, as they are not the main point of this work and as there are no real difficulties in handling them, we leave the proofs to the reader. 

We call \emph{tree} a rooted planar map with only one face. Beware that in particular, the trees we consider are embedded in the plane and that the vertex tree (the only tree with no edges) is allowed here. For $n\ge 0$ and $p\ge 1$, we call \emph{forest} with~$p$ trees and~$n$ edges a $p$-uple of trees where the total number of edges is~$n$. We write $\F_{n,p}$ the set of forests with~$p$ trees and~$n$ edges. A simple application of the so-called cycle lemma (see for example~\cite[Section~6.1]{pitman75} or \cite[Lemma~3]{bettinelli10slr}) yields that 
\begin{equation*}\label{fnp}
|\F_{n,p}|=\frac p {2n+p}\choo{2n+p}{n}.
\end{equation*}

It will be convenient to add extra half-edges between the roots of successive trees: if $\f=(\tr_1,\dots,\tr_p)$ is a forest, we add, for $1\le i \le p$, a half-edge from~$\rho_i$ to~$\rho_{i+1}$, where~$\rho_i$ is the root of~$\tr_i$ and $\rho_{p+1}\de \rho_1$. With this convention, a forest corresponds to a map with two faces, one of which being of degree~$p$ and having a simple curve as a boundary (see Figure~\ref{fnn}). The \emph{corners} of the forest are defined as the corners of the other face of this map. It is easy to see that an element of $\F_{n,p}$ has
\begin{itemize}
	\item $n$ edges,
	\item $n+p$ vertices,
	\item $2n+p$ corners.
\end{itemize}

Before we begin, we need to introduce some vocabulary. For a half-edge~$h$, we write~$h^-$ its origin, $h^+$ its end, and~$\rev(h)$ its reverse. It will be convenient to consider corners as half-edges having no origin, only an end. In particular, if~$c$ is a corner, we will write~$c^+$ the vertex corresponding to it, that is, if~$c$ is the corner delimited by the consecutive half-edges~$h$ and~$h'$, then $c^+\de h^+=h'^-$.

\begin{defi}
A \emph{path} from a vertex~$v$ to a vertex~$v'$ is a finite sequence $\pp=(\pp_1,\pp_2,\dots,\pp_\ell)$ of half-edges such that $\pp_1^-=v$, for $1\le k \le \ell-1$, $\pp_k^+=\pp_{k+1}^-$, and $\pp_\ell^+=v'$. Its \emph{length} is the integer $[\pp]\de\ell$. We will use the convention that an empty path has length~$0$.

A path~$\pp$ is called \emph{self-avoiding} if it does not meet twice the same vertex, that is, $$\big|\{\pp_1^-,\dots,\pp_{[\pp]}^-,\pp_{[\pp]}^+\}\big|=[\pp]+1.$$

The \emph{reverse} of a path~$\pp=(\pp_1,\pp_2,\dots,\pp_\ell)$ is the path $\rev(\pp) \de (\rev(\pp_\ell),\rev(\pp_{\ell-1}),\dots,\rev(\pp_1))$. 
\end{defi}

Let~$\pp$ be a path. We denote by $\pp_{i\to j}$ the path $(\pp_i,\dots,\pp_j)$ if $1\le i\le j\le [\pp]$, or the empty path otherwise. If~$\mathbbm q$ is another path satisfying $\mathbbm q_1^-=\pp_{[\pp]}^+$, we set
$$\pp \ooo \mathbbm q \de (\pp_1,\dots,\pp_{[\pp]},\mathbbm q_1,\dots,\mathbbm q_{[\mathbbm q]})$$
the concatenation of~$\pp$ and~$\mathbbm q$. Throughout this paper, the notion of metric we use is the graph metric: if~$\m$ is a map, the distance $d_\m(v,v')$ between two vertices~$v$ and~$v'$ is the smaller~$\ell$ for which there exists a path of length~$\ell$ from~$v$ to~$v'$. A \emph{geodesic} from~$v$ to~$v'$ is such a path. 

\begin{figure}[ht]
		\psfrag{c}[][]{$c$}
		\psfrag{d}[][]{$c'$}
		\psfrag{o}[][]{or}
	\centering\includegraphics{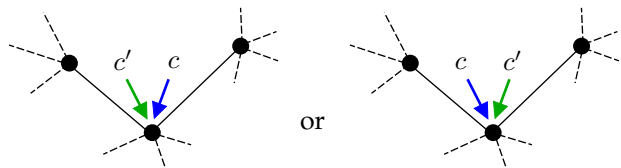}
	\caption{The two different ways of distinguishing twice the same corner.}
	\label{discor}
\end{figure}

We will also say that a half-edge \emph{is directed toward} a set if its end is strictly closer to the set than its origin. A half-edge \emph{is directed away} from a set if its reverse is directed toward the set. In what follows, we will always use the convention that distinguishing a corner ``splits'' it into two new corners. In other words, when we distinguish the same corner for the second time, we have to specify which of its two sides is distinguished (see Figure~\ref{discor}).

\subsection{Adding an edge}\label{secfnn}

Our first bijection is between the set $\F_{n,p}^{\textup{cc}}$ of forests from $\F_{n,p}$ carrying two distinguished corners and the set $\F_{n+1,p}^{\textup{ev}}$ of forests from $\F_{n+1,p}$ having a distinguished edge and a distinguished vertex. It provides a combinatorial interpretation to the following identity:
$$(2n+p)(2n+p+1)\, |\F_{n,p}|=(n+1)(n+p+1)\, |\F_{n+1,p}|.$$

\begin{figure}[ht]
		\psfrag{v}[][][.8]{$v$}
		\psfrag{e}[][][.8]{$\vec e$}
		\psfrag{c}[][][.8]{$c$}
		\psfrag{d}[][][.8]{$c'$}
		\psfrag{l}[][][.8]{\textcolor{red}{$\pl$}}
		\psfrag{r}[][][.8]{\textcolor{red}{$\pr$}}
	\centering\includegraphics[width=.95\linewidth]{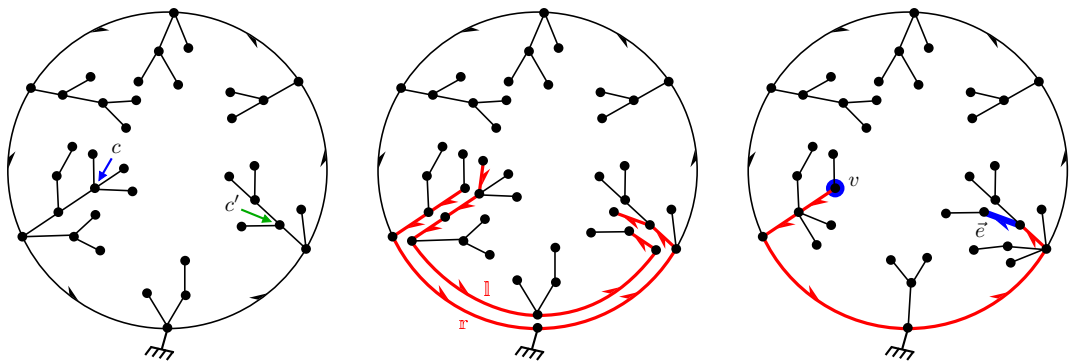}
	\caption{Adding an edge to a forest. The rake-like symbol represents the first tree of the forest.}
	\label{fnn}
\end{figure}

Let $(\f;c,c') \in \F_{n,p}^{\textup{cc}}$. We let~$\pp$ be the only self-avoiding path from~$c^+$ to~$c'^+$ (recall that we added half-edges between the roots of the trees). Starting from the first tree of~$\f$ and following the contour of its trees one by one gives a natural ordering of its corners. A classical way of picturing this order is to imagine an animal flying around the first tree, then jumping to the second one and so on. We ``cut'' along~$\pp$ between~$c$ and~$c'$ in the following sense: we consider, on the one hand, the part of~$\f$ made of all the elements encountered between~$c$ and~$c'$ in its contour and, on the other hand, the part made of all the elements encountered between~$c'$ and~$c$, together with the path~$\pp$. We call these parts respectively \emph{left part} and \emph{right part}. The path~$\pp$ has a copy in each part: in the left part, we consider~$c$ as an extra half-edge and set $\pl\de c \ooo \pp$; in the right part, we consider~$c'$ as an extra half-edge and set $\pr\de \pp \ooo \rev(c')$. Then, we define the forest~$\f'$ by gluing back together both parts while matching~$\pl_k$ with $\pr_k$ for $1 \le k \le [\pl]=[\pr]$. In~$\f'$, we define~$e$ as the edge corresponding to~$\pl_{[\pl]}$ and $v\de \pl_{1}^-$. We set $\Psi_{n\up,p}(\f;c,c')\de (\f';e,v)$. See Figure~\ref{fnn}.

\bigskip

Conversely, let $(\f';e,v) \in \F_{n+1,p}^{\textup{ev}}$. We denote by~$\vec e$ the half-edge corresponding to~$e$ directed away from~$v$, and~$\pp$ the only self-avoiding path from~$v$ to~$\vec e\,^-$. We will cut along $\pp \ooo \vec e$. To this end, we need to specify between which corners we cut. We cut along $\pp \ooo \vec e$ between the corner delimited by~$(\pp \ooo \vec e)_1$ and its predecessor in the contour and the corner delimited by~$\rev(\vec e)$ and its predecessor. In other words, we cut along $\pp \ooo \vec e$, starting from the first corner to the left before $\pp \ooo \vec e$ and stopping at the first corner to the right after $\pp \ooo \vec e$. As above, this defines a left part and a right part. We denote by~$\pl$ the image of $\pp \ooo \vec e$ in the left part and~$\pr$ its image in the right part. We define the forest~$\f$ by gluing back these parts while matching~$\pl_{k+1}$ with~$\pr_k$ for $1\le k \le [\pl]-1$. Because of how we chose the corners between which we cut, in the left part, nothing except~$\pl_1$ is attached to~$\pl_1^-$. As a result, $\pl_1$, which is not glued with anything, gives birth to a corner, which we denote by~$c$. Similarly, $\rev(\pr_{[\pr]})$ defines a corner~$c'$ and we set $\Psi_{n+1\dn,p}(\f';e,v)\de (\f;c,c')$.

\begin{thm}\label{thmfnn}
The mappings $\Psi_{n\up,p} : \F_{n,p}^{\textup{cc}} \to \F_{n+1,p}^{\textup{ev}}$ and $\Psi_{n+1\dn,p} : \F_{n+1,p}^{\textup{ev}} \to \F_{n,p}^{\textup{cc}}$ are one-to-one and reverse one from another.
\end{thm}

\subsection{Adding a tree}

Our second bijection on forests is between the set $\F_{n,p}^{\textup{ci}}$ of forests from $\F_{n,p}$ carrying a distinguished corner together with an integer in $\{1,\dots, p+1\}$ and the set $\F_{n,p+1}^{\textup{vi}}$ of forests from $\F_{n,p+1}$ carrying a distinguished vertex together with an integer in $\{1, \dots,p\}$. It accounts for the following combinatorial identity:
$$(2n+p)(p+1)\, |\F_{n,p}|=(n+p+1)\,p\, |\F_{n,p+1}|.$$

\begin{figure}[ht]
		\psfrag{v}[][][.8]{$v$}
		\psfrag{c}[][][.8]{$c$}
		\psfrag{d}[][][.8]{$c'$}
		\psfrag{l}[][][.8]{\textcolor{red}{$\pl$}}
		\psfrag{r}[][][.8]{\textcolor{red}{$\pr$}}
		\psfrag{1}[][][.8]{$1$}
		\psfrag{2}[][][.8]{$2$}
		\psfrag{3}[][][.8]{$3$}
		\psfrag{4}[][][.8]{$4$}
		\psfrag{5}[][][.8]{$5$}
		\psfrag{6}[][][.8]{$6$}
		\psfrag{7}[][][.8]{$7$}
	\centering\includegraphics[width=.95\linewidth]{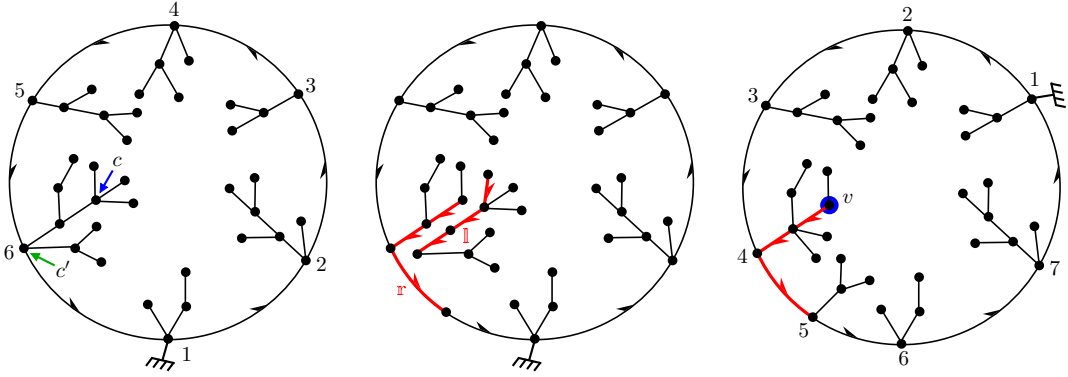}
	\caption{Adding a tree to forest. In this example, $i=4$ and $j=6$.}
	\label{fpp}
\end{figure}

Let us take $(\f;c,i) \in \F_{n,p}^{\textup{ci}}$. The vertex~$c^+$ belongs to some tree of~$\f$: let $j\in \{1,\dots, p\}$ be its index and~$c'$ its last corner. With $(\f;c,c')$, we define the same path~$\pp$ as at the beginning of Section~\ref{secfnn} and we perform the same operation to define a left and a right part. In the left part, we still consider~$c$ as an extra half-edge and set $\pl\de c \ooo \pp$. In the right part, we add an extra tree consisting in one vertex between the $j$-th and the $j+1$-th tree of~$\f$. We denote by~$h$ the half-edge linking the root of the $j$-th tree to this extra tree and define $\pr\de \pp \ooo h$. Then, as above, we glue back together the left part and the right part while matching~$\pl_k$ with $\pr_k$ for $1 \le k \le [\pl]=[\pr]$ and we define $v\de \pl_{1}^-$. Finally, we define~$\f'$ by re-rooting this forest in such a way that~$v$ belongs to the~$i$-th tree and we set $\Psi_{n,p\up}(\f;c,i)\de (\f';v,j)$. See Figure~\ref{fpp}.

Note that the somehow surprising re-rooting accounts for the factors $p+1$ in the left-hand side and~$p$ in the right-hand side of the combinatorial identity. Without this re-rooting, we would miss the forests carrying the vertex~$v$ in their last tree.

Conversely, let us take $(\f';v,j) \in \F_{n,p+1}^{\textup{vi}}$. Let~$i$ be the index of the tree to which~$v$ belongs and consider the path~$\pp$ from~$v$ to the root of the $i+1$-th tree (with the convention that the $i+1$-th tree is the first tree is $i=p+1$). We cut along this path between the corner delimited by~$\pp_1$ and its predecessor in the contour and the last corner of the $i+1$-th tree. This again defines a left part and a right part; we denote respectively by~$\pl$ and~$\pr$ the images of~$\pp$ in these parts. We define a new forest by gluing them back while matching~$\pl_{k+1}$ with~$\pr_k$ for $1\le k \le [\pl]-1$. As before, $\pl_1$ gives birth to a corner~$c$. Moreover, the $i+1$-th tree of this new forest consists in only one vertex; we remove it. Finally, we define~$\f$ by re-rooting the latter forest in such a way that~$c^+$ belongs to the~$j$-th tree, and we set $\Psi_{n,p+1\dn}(\f';v,j)\de (\f;c,i)$.

\begin{thm}\label{thmfpp}
The mappings $\Psi_{n,p\up} : \F_{n,p}^{\textup{ci}} \to \F_{n,p+1}^{\textup{vi}}$ and $\Psi_{n,p+1\dn} : \F_{n,p+1}^{\textup{vi}} \to \F_{n,p}^{\textup{ci}}$ are one-to-one and reverse one from another.
\end{thm}

\section{Modifying the length of the boundary in a quadrangulation}\label{secpp}

We leave forests and concentrate on quadrangulations again. In this section, we will use the corners of the external face. Recall the convention (stated right before Section~\ref{secfnn}) we use for distinguishing corners. As the external face of a quadrangulation from $\Qnp$ already has the root corner distinguished, we will consider that it has $2p+1$ corners. If we distinguish one of these $2p+1$ corners, it will then have $2p+2$ corners. On the one hand, we consider the set $\Qnp^{\textup{ecc}}$ of quadrangulations from~$\Qnp$ carrying
\begin{itemize}
	\item one distinguished edge~$e$,
	\item one first distinguished corner~$c$ of the external face,
	\item one second distinguished corner~$c'$ of the external face.
\end{itemize}
On the other hand, we consider the set $\Q_{n,p+1}^{\textup{vhh}}$ of quadrangulations from~$\Q_{n,p+1}$ carrying
\begin{itemize}
	\item one distinguished vertex~$v$,
	\item one first half-edge~$h$ of the boundary directed toward~$v$,
	\item one second half-edge~$h'\neq h$ of the boundary directed toward~$v$.
\end{itemize}
Note that there are exactly $p+1$ half-edges of the boundary that are directed toward~$v$. To see this, label the vertices of the boundary with their distance to~$v$. As there are no cycles of odd length in a quadrangulation, the labels encountered when traveling along the boundary form a $2(p+1)$-step bridge whose steps are either~$+1$ or~$-1$. As a result, exactly half of them are~$-1$ steps.

We will present an explicit bijection between the previous two sets; this will provide a combinatorial interpretation to~\eqref{pp}.

\bigskip

We will often use the notion of \emph{left-most geodesic} from a half-edge~$h$ (or a corner) to some set~$S$ of vertices (typically a vertex, an edge or a face). It is constructed as follows. First, we consider all the geodesics from~$h^+$ to the closest elements of~$S$. We take the set of all the first steps of these geodesics. Starting from~$h$, we select the first half-edge to its left that belongs to this set. In other words, we turn clockwise around $h^+$ and select the first half-edge of this set that we meet. Note that this half-edge may be $\rev(h)$ if this is the only half-edge in the set. Then we iterate the process from this half-edge until we reach~$S$. Note that this path may be empty if~$h^+\in S$ and that it is a geodesic. The \emph{right-most geodesic} from a half-edge or a corner to some set is defined in a similar way, by changing left to right in the previous definition.

\subsection{Adding two half-edges to the boundary}\label{secpp+}

\begin{figure}[ht]
		\psfrag{c}[][]{$c$}
		\psfrag{e}[][]{$e$}
	\centering\includegraphics{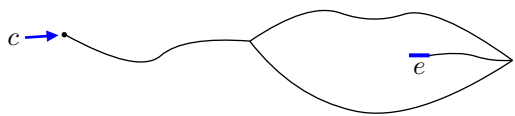}
	\caption{Example where the right-most geodesic from~$\vec e$ to~$c$ (on the bottom) is not the reverse of the left-most geodesic from~$c$ to~$e$ (on the top).}
	\label{revc}
\end{figure}

We begin with the mapping from $\Qnp^{\textup{ecc}}$ to $\Q_{n,p+1}^{\textup{vhh}}$. Let $(\q;e,c,c')\in \Qnp^{\textup{ecc}}$. Note that, as~$\q$ is a quadrangulation, one of the extremities of~$e$ is strictly closer to~$c^+$ than the other. We denote by~$\vec e $ the half-edge corresponding to~$e$ directed toward~$c^+$, and~$\pc$ the right-most geodesic from~$\vec e$ to~$c^+$. Beware that $\rev(\pc)$ is not necessarily the left-most geodesic from~$c$ to~$e$: there may exist a path of same length leaving~$\pc$ to the left and meeting it again from the right, as shown on Figure~\ref{revc}.

\begin{figure}[ht]
		\psfrag{v}[][][.8]{$v$}
		\psfrag{h}[][][.8]{$h$}
		\psfrag{i}[][][.8]{$h'$}
		\psfrag{a}[][][.8]{\textcolor{red}{$\pl$}}
		\psfrag{b}[][][.8]{\textcolor{red}{$\pr$}}
		\psfrag{c}[][]{$c$}
		\psfrag{d}[][]{$c'$}
		\psfrag{x}[][][.8]{\textcolor{red}{$\ph_0$}}
		\psfrag{y}[l][l][.8]{\textcolor{red}{$\ph'_0$}}
		\psfrag{e}[][]{\textcolor{blue}{$e$}}
		\psfrag{1}[][][.85]{$\pp_1$}
		\psfrag{2}[][][.85]{$\pp_2$}
		\psfrag{3}[][][.85]{$\pp_3$}
		\psfrag{4}[][][.85]{$\pp_4$}
		\psfrag{5}[][][.85]{$\pp_5$}
		\psfrag{6}[][][.85]{$\pp_6$}
		\psfrag{7}[][][.85]{$\pp_7$}
		\psfrag{8}[][][.85]{$\pp_8$}
	\centering\includegraphics[width=12cm]{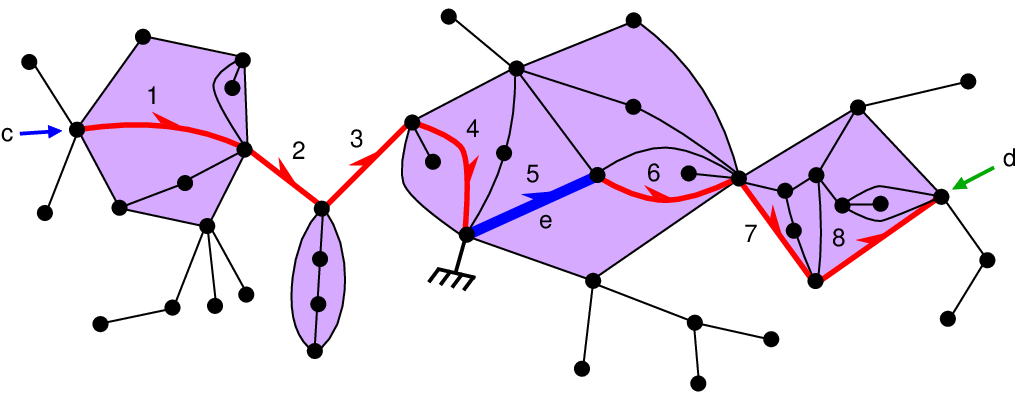}
	\bigskip
	
	\centering\includegraphics[width=.95\linewidth]{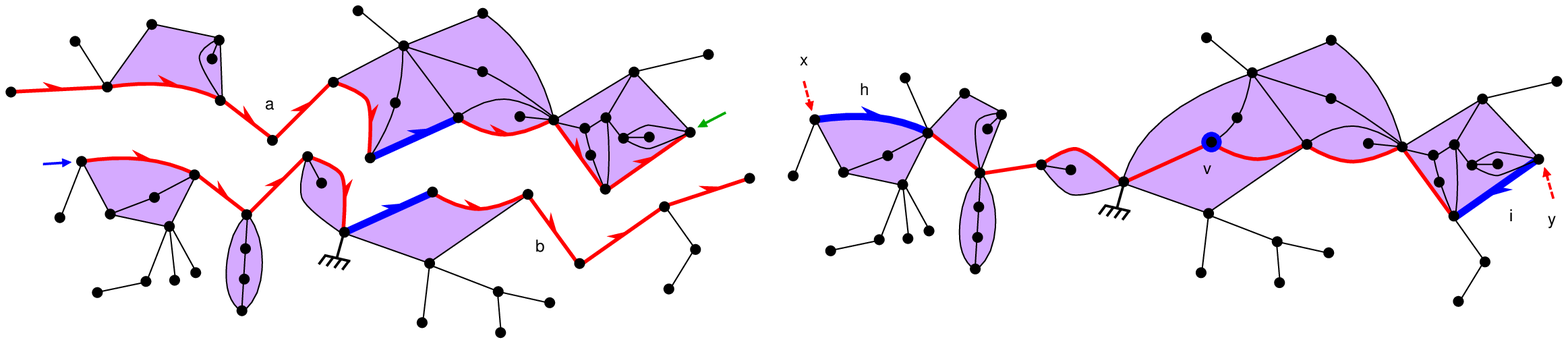}
	\caption{The mapping from $\Qnp^{\textup{ecc}}$ to $\Q_{n,p+1}^{\textup{vhh}}$ in the simple case. We define the path~$\pp$, then cut along it and glue back after shifting the parts one notch. This creates two new half-edges on the boundary, $h$ and~$h'$. On this example, $[\pc]=4$ and $[\pc']=3$.}
	\label{pp+}
\end{figure}

Until further notice, we suppose that $\rev(\vec e)$ is directed toward~$c'^+$: in this case, we say that the quadruple $(\q;e,c,c')$ is \emph{simple}. We denote by~$\pc'$ the right-most geodesic from~$\rev(\vec e)$ to~$c'^+$. It is not hard to see that, in this case, $\rev(\pc)$ is actually the left-most geodesic from~$c$ to~$e$: indeed, a path leaving~$\rev(\pc)$ to the left and meeting it from the right has to intersect~$\pc'$; as a result, this path cannot be a geodesic as $\rev(\vec e)$ is directed toward~$c'^+$. We split the map~$\q$ in two parts between~$c$ and~$c'$ along the self-avoiding path
$$\pp \de \rev(\pc)\ooo \rev(\vec e) \ooo \pc'.$$
Let us call \emph{left part} (resp.\ \emph{right part}) the part consisting of the path~$\pp$ together with all the elements located to its left (resp.\ its right). In the left part, we replace~$c$ with an extra half-edge~$h$ and set $\pl\de h \ooo \pp$. In the right part, we replace~$c'$ with an extra half-edge~$h'$ and set $\pr\de \pp \ooo \rev(h')$. Then, we glue back the two parts together in such a way that $\pl_k$ coincides with~$\pr_k$ for $1\le k \le [\pl]=[\pr]$. We denote by~$\q'$ the map we obtain and we set $v\de\pl_{[\pc]+1}^+ = (\rev(\pr))_{[{\pc'}]+1}^+$. See Figure~\ref{pp+}.

\begin{rem}
Note that, when sliding one notch along the path, the marked edge is duplicated and the intersection between the two copies of this edge will become the new marked vertex. Conversely, a marked vertex will be duplicated into two vertices delimiting a new marked edge. This general principle will often be used in what follows.
\end{rem}

\bigskip

\begin{figure}[ht]
		\psfrag{v}[][][.8]{$v$}
		\psfrag{h}[][][.8]{$h$}
		\psfrag{i}[][][.8]{$h'$}
		\psfrag{a}[][][.8]{\textcolor{red}{$\pl$}}
		\psfrag{b}[][][.8]{\textcolor{red}{$\pr$}}
		\psfrag{c}[][][.8]{$c$}
		\psfrag{d}[][][.8]{$c'$}
		\psfrag{x}[][][.8]{\textcolor{red}{$\ph_0$}}
		\psfrag{y}[l][l][.8]{\textcolor{red}{$\ph'_0$}}
		\psfrag{e}[][][.8]{\textcolor{blue}{$e$}}
		\psfrag{1}[][][.75]{$\pp_1$}
		\psfrag{2}[][][.75]{$\pp_2$}
		\psfrag{3}[][][.75]{$\pp_3$}
		\psfrag{4}[][][.75]{$\pp_4$}
		\psfrag{5}[][][.75]{$\pp_5$}
		\psfrag{6}[][][.75]{$\pp_6$}
		\psfrag{7}[][][.75]{$\pp_7$}		
	\centering\includegraphics[width=.95\linewidth]{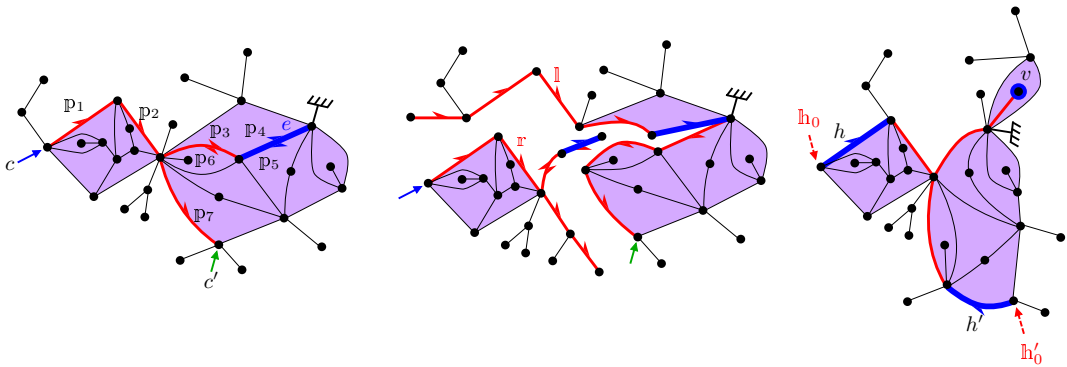}
	\caption{The shifting operation in the left-pinched case. We shift along the path~$\pp$ and circumvent the part where~$\pp$ is pinched. On this picture, $[\pc]=3$ and $[\pc']=2$.}
	\label{pp+2}
\end{figure}

Now, let us suppose that~$\vec e$ is directed toward~$c'^+$. We denote by~$\pc'$ the right-most geodesic from~$\vec e$ to~$c'^+$. Let $i\ge 1$ be the smallest integer such that $\pc_i \neq \pc'_i$. Note that~$\pc$ and~$\pc'$ do not meet again after time~$i$, in the sense that $\{\pc_i^+,\dots,\pc_{[\pc]}^+\}$ and $\{\pc_i'^+,\dots, \pc_{[\pc']}'^+\}$ are disjoint sets. If~$\pc_i'$ is to the left of~$\pc_i$ after~$\pc_{i-1}$ (with the convention $\pc_0=\pc_0' \de \vec e$), we say that the quadruple $(\q;e,c,c')$ is \emph{left-pinched}; we say that it is \emph{right-pinched} in the other case. This terminology comes from the fact that the path
$$\pp \de \rev(\pc) \ooo \rev(\vec e) \ooo \vec e \ooo \pc'$$
is ``pinched,'' the pinched part being to the left or right of the path $\rev(\pc_{i\to [\pc]}) \ooo \pc'_{i \to [\pc']}$. Let us suppose here that $(\q;e,c,c')$ is left-pinched. Similarly as before, we split the map in two parts. The right part is defined as previously. The definition of the left part, however, is slightly different: the difference is that we cut along the pinched part $\vec e \ooo \pc_{1\to i-1}$ (see the middle part of Figure~\ref{pp+2}). We define~$h$, $h'$ and the paths~$\pl$ and~$\pr$ by the same method as above. Note that here~$\pl$ is self-avoiding whereas~$\pr$ is not. We then proceed as before: we glue the parts while matching $\pl_{k}$ with~$\pr_k$ for $1\le k \le [\pl]$. Let~$\q'$ be the map we obtain and let us define $v\de\pl_{[\pc]+1}^+ = (\rev(\pr))_{[{\pc'}]+2}^+$. See Figure~\ref{pp+2}.

To treat the right-pinched case, we use a symmetric argument, by exchanging the roles of~$c$ and~$c'$. Note that $(\q;e,c,c')$ is right-pinched if and only if $(\q;e,c',c)$ is left-pinched. We define the quadruple $(\q';v,h',h)$ corresponding to $(\q;e,c',c)$ by the previous construction. We then associate with $(\q;e,c,c')$ the quadruple $(\q';v,h,h')$. In all cases, we set $\Phi_{n,p\up}(\q;e,c,c') \de (\q';v,h,h')$.

\begin{rem}
As a forest may be seen as a particular $2$-face map, one might wonder whether $\Phi_{n,p\up}$ is a simple extension of the mapping presented in Section~\ref{secfnn}. This is not the case as its restriction $\Phi_{0,p\up}$ to trees is not the same mapping as $\Psi_{p\up,1}$ from Section~\ref{secfnn}. Interestingly, this provides yet another mapping on trees, leaving us with three different mappings: R\'emy's algorithm, $\Phi_{0,p\up}$ and $\Psi_{p\up,1}$.
\end{rem}

\subsection{Removing two half-edges from the boundary}\label{secpp-}

Let us now describe the reverse mapping $\Phi_{n,p+1\dn}$. We will also have to distinguish between three cases, and we try to keep our notation coherent with the previous section. Let $(\q';v,h,h') \in \Q_{n,p+1}^{\textup{vhh}}$. Let~$\ph_0$ denote the corner delimited by~$h$ and its predecessor in the contour of the external face, and~$\ph$ the left-most geodesic from this corner to~$v$. Note that, as~$h$ is directed toward~$v$ and is the first half-edge to the left after~$\ph_0$, we necessarily have $[\ph]\ge1$ and $\ph_1=h$. We define~$\ph'_0$ and~$\ph'$ in a similar way with~$h'$ instead of~$h$. As $h\neq h'$, there is no ambiguity in the definition of~$\ph_0$ and~$\ph'_0$.

If~$\ph$ and~$\ph'$ do not meet before reaching~$v$, 
we say that $(\q';v,h,h')$ is \emph{simple}. In this case, we set $\pp \de \ph \ooo \rev(\ph')$. Notice that this path is self-avoiding. As before, we define a left part and a right part by cutting along~$\pp$ between the corners~$\ph_0$ and~$\ph'_0$. We denote by~$\pl$ and~$\pr$ the path~$\pp$ in these respective parts. We define~$\q$ by gluing back the two parts while matching~$\pl_{k+1}$ with~$\pr_k$ for $1\le k \le [\pl]-1$. As in Section~\ref{secfnn}, $\pl_1$ and $(\rev(\pr))_1$ define two corners~$c$ and~$c'$ in~$\q$. Finally, we let~$e$ be the edge corresponding to~$\pl_{[\ph]+1}$ and~$(\rev(\pr))_{[\ph']+1}$.

\bigskip

If~$\ph$ and~$\ph'$ meet before reaching~$v$, we let~$i$ and~$j$ be the smallest integers such that $\ph_{i}^+=\ph_{j}'^+$. Note that~$\ph$ and~$\ph'$ merge after times~$i$ and~$j$, in the sense that $\ph_{i+k} = \ph_{j+k}'$ for all $k\ge 1$ such that these quantities are defined. To see this, observe that the ``path'' $\ph_0\ooo\ph_{1\to i}\ooo \rev(\ph'_{1\to j})\ooo\rev(\ph'_0)$ separates the map into two disjoint components and that~$v$ belongs to only one of them. We use the same vocabulary as in the previous section, left-pinched or right-pinched, depending on whether~$v$ belongs to the component located to the left or to the right of this path.

Let us suppose that $(\q';v,h,h')$ is left-pinched. In this case, it is not hard to see that $j\ge1$
. We define $\pp \de \ph \ooo \rev(\ph')$ and, as in the previous section, a left part with a self-avoiding path~$\pl$ and a right part with a path~$\pr$. The map~$\q$ is then obtained by gluing $\pl_{k+1}$ to $\pr_{k}$ for $1\le k \le [\pl]-1$. The fact that $j\ge 1$ shows that only $(\rev(\pr))_1$ is attached to~$(\rev(\pr))_1^-$ in the right part, so that~$(\rev(\pr))_1$ creates a corner~$c'$. As~$\pl$ is self-avoiding, $\pl_1$ also defines a corner~$c$. The edge~$e$ is the one corresponding to~$\pr_{[\ph]}$.

The right-pinched case is treated by the same method as in the previous section, by exchanging the roles of~$h$ and~$h'$. We define $\Phi_{n,p+1\dn}(\q';v,h,h')\de(\q;e,c,c')$.

\subsection{These mappings are reverse one from another}

Now that we have described in detail our mappings, we may show the main result of this section:

\begin{thm}\label{thmpp}
The mappings $\Phi_{n,p\up}: \Qnp^{\textup{ecc}} \to \Q_{n,p+1}^{\textup{vhh}}$ and $\Phi_{n,p+1\dn}: \Q_{n,p+1}^{\textup{vhh}} \to \Qnp^{\textup{ecc}}$ are one-to-one and reverse one from another.
\end{thm}

\begin{pre}
We will see that simple quadruples of $\Qnp^{\textup{ecc}}$ correspond to simple quadruples of $\Q_{n,p+1}^{\textup{vhh}}$ through both mappings, and that the same goes for left-pinched quadruples and right-pinched quadruples. We treat these cases separately.

\paragraph{Simple case.} 
Let $(\q;e,c,c')\in\Qnp^{\textup{ecc}}$ be a simple quadruple and $(\q';v,h,h') \de \Phi_{n,p\up}(\q;e,c,c')$. We claim that it is a simple quadruple of $\Q_{n,p+1}^{\textup{vhh}}$. The fact that $\q'\in \Q_{n,p+1}$ is clear, as $\Phi_{n,p\up}$ only alters the external face by adding the two half-edges~$h$ and~$h'$. For now, we use the notation of Section~\ref{secpp+}. We claim that, in~$\q'$, the path $\pl_{1\to [\pc]+1}$ is the left-most geodesic from the last corner before~$h$ in the contour toward~$v$. To show this, we use an argument that we will often use throughout this paper. We argue by contradiction and suppose that this does not hold. Let~$\ph$ (with possibly $[\ph] < [\pc]+1$) be the left-most geodesic in question. We will see that~$\ph$ has to ``wind'' an infinite number of times around~$v$. Let $k\le [\pc]+1$ be the smallest integer such that $\ph_k\neq \pl_k$. Let us first suppose that~$\ph_k$ belongs to the left part. The path~$\ph$ cannot stay in the left part as otherwise it would hit the set $\{\pl_1^+,\dots,\pl_{[\pc]+1}^+\}$ (to which~$v$ belongs) and, in~$\q$, this would either provide a better alternative to~$\pc$ or contradict the fact that~$\vec e$ is directed toward~$c^+$ (this case could happen if~$\ph$ ends with $\rev(\pl_{[\pc]+2})$). Let~$l>k$ be the smallest integer such that~$\ph_l$ does not lie in the left part. By the preceding argument, we see that $\ph_l^- \in \{\pl_{[\pc]+2}^+ , \dots, \pl_{[\pl]}^+\}$. A symmetric argument shows that~$\ph$ has to leave the right part: we let~$m$ be the smallest integer greater than~$l$ such that~$\ph_m$ does not lie in the right part. Then, we have $\ph_m^-\in\{\pl_k^+,\dots,\pl_{[\pc]}^+\}$. As a result, $\ph_{k\to m-1}$ creates a wind around~$v$ and we are back to the same situation as before with~$\ph_m$ instead of~$\ph_k$. Reiterating the argument, we obtain that~$\ph$ has to infinitely wind around~$v$, which is a contradiction. By similar arguments, we obtain that~$\ph_k$ does not lie in the right part either and our claim follows.

\begin{figure}[ht]
		\psfrag{k}[][]{$\ph_k$}
		\psfrag{u}[][]{$\ph_l$}
		\psfrag{m}[][]{$\ph_m$}
		\psfrag{v}[][]{$v$}
		\psfrag{h}[][]{$h$}
		\psfrag{i}[][]{$h'$}
		\psfrag{l}[][]{\textcolor{red}{$\pl$}}
	\centering\includegraphics{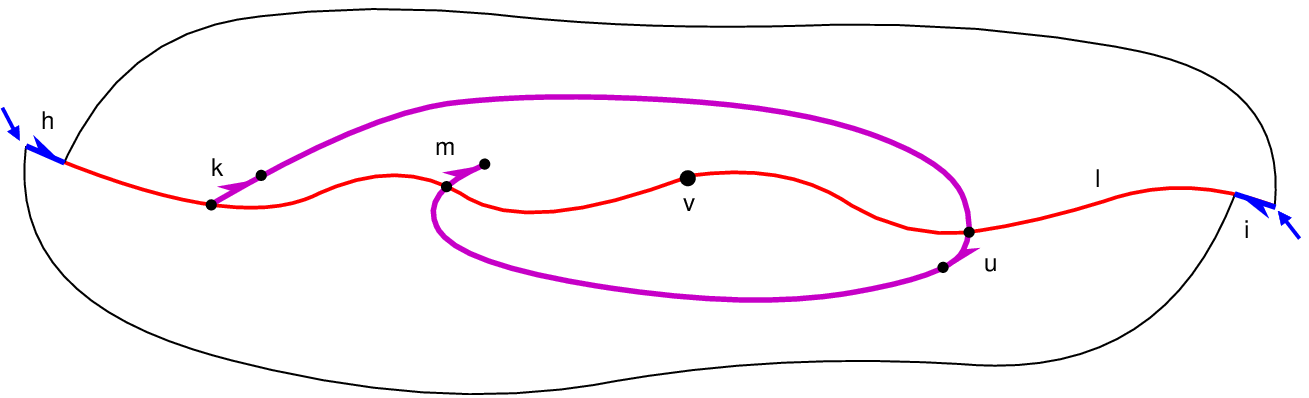}
	\caption{The first wind of~$\ph$ around~$v$.}
	\label{loop}
\end{figure}

By symmetry, the path $\rev(\pr)_{1\to [\pc']+1}$ is the left-most geodesic from the last corner before~$h'$ to~$v$. In particular, $h$ and~$h'$ are directed toward~$v$, so that $(\q';v,h,h')\in \Q_{n,p+1}^{\textup{vhh}}$. With the notation of Section~\ref{secpp-}, we see that these two paths are~$\ph$ and~$\ph'$, so that the notation~$\pl$ and~$\pr$ from both sections is coherent. It becomes clear that $(\q';v,h,h')$ is simple and that $\Phi_{n,p+1\dn}(\q';v,h,h') =(\q;e,c,c')$.

Conversely, let $(\q';v,h,h') \in \Q_{n,p+1}^{\textup{vhh}}$ be simple, and $(\q;e,c,c')\de \Phi_{n,p+1\dn} (\q';v,h,h')$. A similar method shows that, with the notation of Section~\ref{secpp-}, the path $\rev(\pl_{2 \to [\ph]})$ is the right-most geodesic in~$\q$ from~$\vec e$ to~$c^+$. Using a symmetry argument, this is sufficient to conclude that $(\q;e,c,c')$ is simple and that $(\q';v,h,h') = \Phi_{n,p\up}(\q;e,c,c')$.

\paragraph{Pinched case.}
As we pass from a left-pinched quadruple to a right-pinched quadruple by exchanging the last two coordinates in both $\Qnp^{\textup{ecc}}$ and $\Q_{n,p+1}^{\textup{vhh}}$, we may restrict our attention to the left-pinched case. Let $(\q;e,c,c')\in\Qnp^{\textup{ecc}}$ be left-pinched and $(\q';v,h,h') \de \Phi_{n,p\up}(\q;e,c,c')$. Let~$\ph_0$ and~$\ph'_0$ be the corners before~$h$ and~$h'$ in the contour of the external face. To conclude that $(\q';v,h,h')$ belongs to $\Q_{n,p+1}^{\textup{vhh}}$, is left-pinched and that $\Phi_{n,p+1\dn}(\q';v,h,h') = (\q;e,c,c')$, it is sufficient to show that, with the notation of Section~\ref{secpp+}, $\pl_{1\to [\pc]+1}$ and $\rev(\pr)_{1\to [\pc']+2}$ are the left-most geodesics from~$\ph_0$ and~$\ph'_0$ to~$v$ in~$\q'$.

We use the same winding argument as in the simple case but some extra care is needed in this case. See Figure~\ref{wind}. Let~$\ph$ be the left-most geodesic from~$\ph_0$ to~$v$. We argue by contradiction and suppose that $\ph \neq\pl_{1\to [\pc]+1}$. Let $k\le [\pc]+1$ be the smallest integer such that $\ph_k\neq \pl_k$. Let us first suppose that~$\ph_k$ belongs to the left part. There are several different ways in which~$\ph$ may hit~$\pl$ after time~$k$: in~$\q'$, let us define
\begin{align*}
i &\de \inf\big\{r: \ \pl_r^+\in \{\pl_1^+,\dots,\pl_{r-1}^+,\pl_{r+1}^+,\dots,\pl^+_{[\pl]}\} \big\} \quad\text{ and}\\
j &\de \sup\big\{r: \ \pl_r^+\in \{\pl_1^+,\dots,\pl_{r-1}^+,\pl_{r+1}^+,\dots,\pl^+_{[\pl]}\} \big\}.
\end{align*}
These are the smallest and largest integers such that $\pl_i^+=\pl_j^+$ in~$\q'$; they separate the path~$\pl$ in four parts, $\pl_{1\to i}$, $\pl_{i+1\to [\pc]+1}$, $\pl_{[\pc]+2 \to j}$ and  $\pl_{j+1 \to [\pl]}$, the second and third being reverse one from another. Beware that in~$\q$, this is no longer true. We consider the smallest $l\ge k$ such that either~$\ph_l$ does not lie in the left part or $\ph_l^+ \in \{\pl_1^+,\dots,\pl_{j}^+\}$. Such an integer exists as~$v$ lies in the latter set. We claim that~$\ph_l$ does not lie in the left part. Indeed, let us argue by contradiction and suppose that $\ph_{k\to l}$ entirely lies in the left part. Let~$s$ be such that $\ph_l^+=\pl_s^+$ and $\pl_s$, $\ph_l$, $\pl_{s+1}$ are arranged according to the clockwise order around $\ph_l^+$. By the following arguments, we obtain a contradiction.
\begin{itemize}
	\item If $s\le [\pc]+1$, in~$\q$, the path $\rev(\ph_{k\to l})$ is to the right of $\rev(\pl_{k\to s})$ and is of the same length or shorter, so that~$\pc$ is not the right-most geodesic from~$\vec e$ to~$c^+$.
	\item If $s=[\pc]+2$, then $[\ph_{k\to l}] \le [\pl_{k\to s}]-2$. As a result, $\vec e$ is not directed toward~$c^+$. 
	\item If $[\pc]+3 \le s \le j$, $\ph_{k\to l}$ creates a shortcut in~$\q$ and~$\pc$ is not a geodesic.
\end{itemize}
As a result, we obtain that~$\ph_l$ does not belong to the left part and that $\ph_l^- \in \{\pl_{j+1}^+ , \dots, \pl_{[\pl]}^+\}$. Let~$m$ be the smallest integer greater than~$l$ such that either~$\ph_m$ does not lie in the right part or $\ph_m^+ \in \{\pl_j^+,\dots,\pl_{[\pl]}^+\}$. It is easy to see that~$\ph_m$ does not lie in the right part, otherwise this would contradict the fact that, in~$\q$, $\pc'$ is the right-most geodesic from~$\vec e$ to~$c'^+$. We then conclude as in the simple case that~$\ph$ has to infinitely wind around~$v$, which is a contradiction.

\begin{figure}[ht]
	\psfrag{0}[][]{$\ph_0$}
		\psfrag{1}[][]{$\ph_0'$}
		\psfrag{k}[][]{$\ph_k$}
		\psfrag{u}[][]{$\ph_l$}
		\psfrag{m}[][]{$\ph_m$}
		\psfrag{t}[][]{$\ph_l$}
		\psfrag{v}[][]{$v$}
		\psfrag{h}[][]{$h$}
		\psfrag{i}[][]{$h'$}
		\psfrag{l}[][]{\textcolor{red}{$\pl$}}
		\psfrag{a}[][]{\textcolor{red}{$\pl_i$}}
		\psfrag{b}[][]{\textcolor{red}{$\pl_j$}}
		\psfrag{s}[][]{\textcolor{red}{$\pl_s$}}
	\centering\includegraphics[width=.9\linewidth]{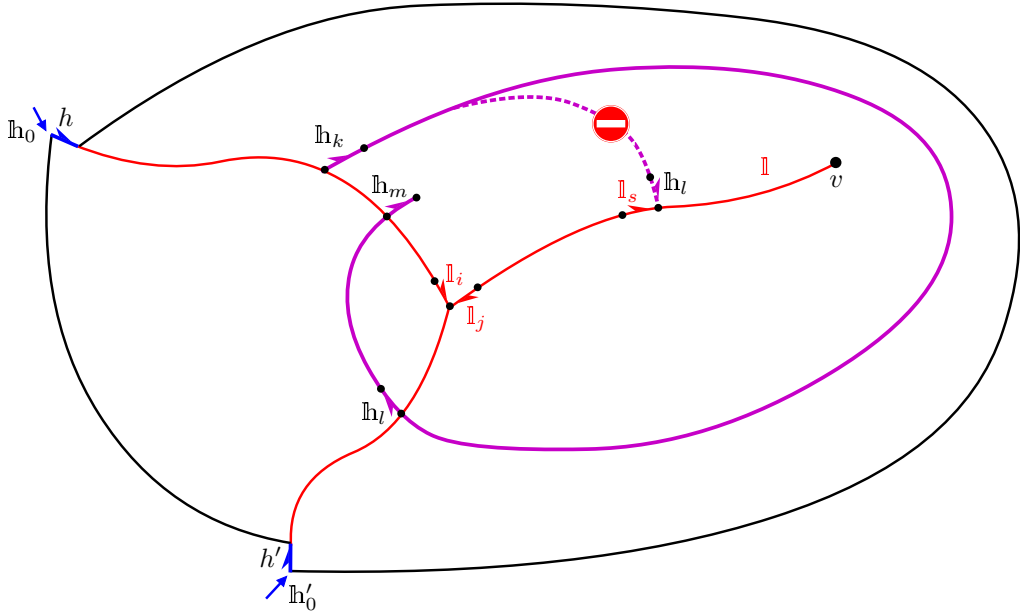}
	\caption{The winding argument in the pinched case.}
	\label{wind}
\end{figure}

By similar arguments, we obtain that~$\ph_k$ does not lie in the right part either, so that $\pl_{1\to [\pc]+1}$ is the left-most geodesic from~$\ph_0$ to~$v$ in~$\q'$. We also obtain that $\rev(\pr)_{1\to [\pc']+2}$ is the left-most geodesic from~$\ph'_0$ to~$v$. This allows us to conclude that $\Phi_{n,p\up}$ maps left-pinched quadruples of $\Qnp^{\textup{ecc}}$ to left-pinched quadruples of $\Q_{n,p+1}^{\textup{vhh}}$ and that $\Phi_{n,p+1\dn}\circ\Phi_{n,p\up}$ is the identity on the set of left-pinched quadruples of $\Qnp^{\textup{ecc}}$. A very similar technique also shows that $\Phi_{n,p+1\dn}$ maps left-pinched quadruples of $\Q_{n,p+1}^{\textup{vhh}}$ to left-pinched quadruples of $\Qnp^{\textup{ecc}}$ and that $\Phi_{n,p\up}\circ\Phi_{n,p+1\dn}$ is the identity on the set of left-pinched quadruples of $\Q_{n,p+1}^{\textup{vhh}}$. We leave the details to the reader.
\end{pre}

\section{Changing the number of faces in a quadrangulation}\label{secnn}

\begin{figure}[ht]
		\psfrag{e}[][]{$e$}
		\psfrag{f}[][]{$e'$}
		\psfrag{o}[][]{or}
	\centering\includegraphics{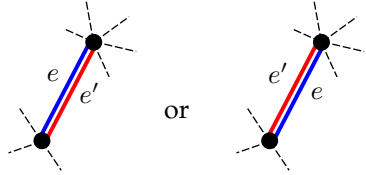}
	\caption{The two different ways of distinguishing twice the same edge.}
	\label{disedge}
\end{figure}

We now present the second main bijection of this work. When distinguishing two edges, we use a convention similar to the one we used for corners. The second time we distinguish an edge, we have to specify on which side it is distinguished (see Figure~\ref{disedge}).
On the one hand, we consider the set $\Qnp^{\textup{eem}}$ of quadruples $(\q;e,e',m)$ where $\q\in \Qnp$ and
\begin{itemize}
	\item $e$ is a distinguished edge of~$\q$,
	\item $e'$ is a second distinguished edge of~$\q$ with the convention mentioned above,
	\item $m$ is a mark in $\{c,s_1,s_2\}$.
\end{itemize}
On the other hand, the set $\Q_{n+1,p}^{\textup{fv}}$ is the set of triples $(\q';f,v)$ where $\q' \in \Q_{n+1,p}$ and
\begin{itemize}
	\item $f$ is a distinguished face of~$\q'$,
	\item $v$ is a distinguished vertex of~$\q'$.
\end{itemize}
Once we will have presentd our bijection between these two sets, we will obtain a proof and combinatorial interpretation of~\eqref{nn}.

\subsection{Removing a face}\label{secnn-}

It seems more pedagogical to start with the mapping $\Phi_{n+1\dn,p}$ from $\Q_{n+1,p}^{\textup{fv}}$ to $\Qnp^{\textup{eem}}$. Let $(\q';f,v) \in \Q_{n+1,p}^{\textup{fv}}$. Using the face~$f$ and the vertex~$v$, we will define a path and shift a part of the map along it in order to suppress~$f$. We will consider two cases: as there are no cycles of odd length in~$\q'$, the distances between the ends of the four corners of~$f$ and~$v$ may only be of two types, either $d$, $d+1$, $d$, $d+1$ or $d$, $d+1$, $d+2$, $d+1$, for some $d\ge 0$. In the first case, we will say that~$f$ is \emph{confluent} with respect to~$v$ and, in the second case, we will say that it is \emph{simple}. This terminology is borrowed from~\cite{schaeffer98cac}.

\subsubsection{Confluent face}\label{secnn-cf}

We suppose here that~$f$ is confluent with respect to~$v$ and we denote by~$d$ the distance between~$f$ and~$v$ (as in the previous definition). Let~$h_1$, $h_2$, $h_3$, $h_4$ be the four half-edges incident to~$f$ read in the counterclockwise order, with~$h_2$ and~$h_4$ directed toward~$v$ (we arbitrarily choose among the two possibilities). We denote by~$\pp$ and~$\pp'$ the left-most geodesics from~$h_2$ and~$h_4$ to~$v$. It is easy to see that~$\pp$ and~$\pp'$ merge as soon as they meet. Note that~$\pp$ and~$\pp'$ may even be equal if~$h_2^+=h_4^+$.

Let us first suppose that~$\pp$ and~$\pp'$ only meet when reaching~$v$. The paths $\rev(h_3) \ooo \pp \ooo \rev(\pp') \ooo \rev(h_4)$ and $h_1 \ooo h_2 \ooo \pp \ooo \rev(\pp')$ are both simple loops, in the sense that they do not go twice through the same vertex\footnote{Note that this remains true if $h_2^+=h_4^+$. In this case, as we supposed that~$\pp$ and~$\pp'$ do not meet before reaching~$v$, we must have $v=h_2^+$.}; as a result, they separate the map in two parts. We call left part the loop $\rev(h_3) \ooo \pp \ooo \rev(\pp') \ooo \rev(h_4)$ together with all the elements of the map located to the left of this loop. In the case where $h_3=\rev(h_4)$, the path $\rev(h_3)\ooo \rev(h_4)$ is a ``flat'' simple loop and the left part is defined as the map consisting in a single edge joining two vertices. In this part, we denote by~$\pl$ the loop. The right part is the loop $h_1 \ooo h_2 \ooo \pp \ooo \rev(\pp')$ together with all the elements of the map located to its right. We use the same convention as above if $h_1=\rev(h_2)$. In the right part, the loop is denoted by~$\pr$. Note that~$f$ belongs to neither parts.

\begin{figure}[ht]
		\psfrag{x}[][][.8]{$h_1$}
		\psfrag{y}[][][.8]{$h_2$}
		\psfrag{v}[][][.8]{$v$}
		\psfrag{f}[][][.8]{$f$}
		\psfrag{m}[][][.8]{$m=c$}
		\psfrag{l}[][][.8]{\textcolor{red}{$\pl$}}
		\psfrag{r}[][][.8]{\textcolor{red}{$\pr$}}
		\psfrag{e}[][][.8]{\textcolor{blue}{$e$}}
		\psfrag{g}[][][.8]{\textcolor{blue}{$e'$}}
		\psfrag{h}[][][.8]{\textcolor{red}{$e'$}}
		\psfrag{1}[][][.75]{$\pl_1$}
		\psfrag{2}[][][.75]{$\pl_2$}
		\psfrag{3}[][][.75]{$\pl_3$}
		\psfrag{4}[][][.75]{$\pl_4$}
		\psfrag{5}[][][.75]{$\pl_5$}
		\psfrag{6}[][][.75]{$\pl_6$}
		\psfrag{7}[][][.75]{$\pr_1$}
		\psfrag{8}[][][.75]{$\pr_2$}
		\psfrag{9}[][][.75]{$\pr_3$}
		\psfrag{a}[][][.75]{$\pr_4$}
		\psfrag{b}[][][.75]{$\pr_5$}
		\psfrag{c}[][][.75]{$\pr_6$}	
	\centering\includegraphics[width=.95\linewidth]{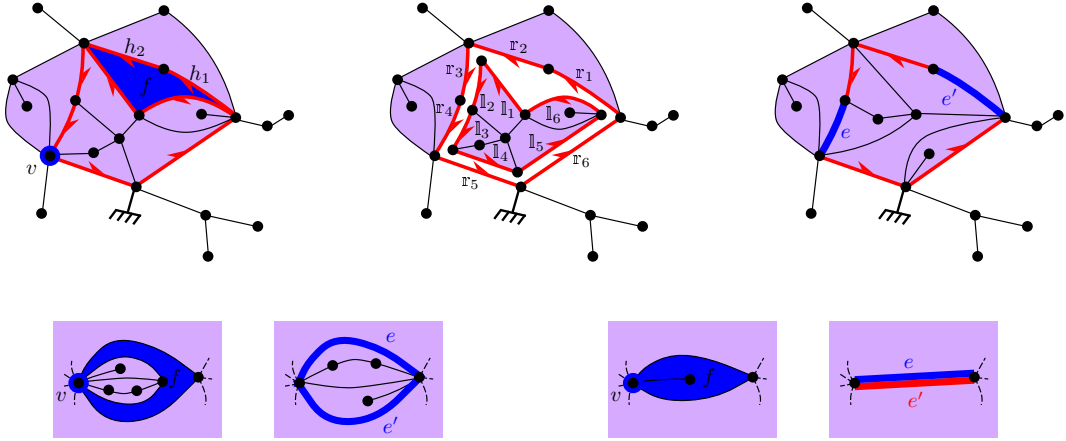}
	\caption{Removing a confluent face. On the first line, $d=2$. The spacial cases where $h_2^+=h_4^+$ and where two half-edges incident to~$f$ are reverse one from another are illustrated on the second line. Beware that, when $h_2^+=h_4^+$ but neither $h_1=\rev(h_2)$ nor $h_3 = \rev(h_4)$, the two parts separated by~$f$ are rotated before being glued back together.}
	\label{nn-}
\end{figure}

The map~$\q$ is the map obtained by gluing back the two parts while matching~$\pl_k$ with~$\pr_k$, for all $1\le k \le [\pl]=2d+2$. The edge~$e$ is the one corresponding to $\pl_{d+2}$ and~$e'$ is the one corresponding to~$\pl_1$. Remark that when $h_2^+=h_4^+$, $e$ and~$e'$ have the same extremities. When, moreover, two half-edges incident to~$f$ are reverse one from another, namely $h_1=\rev(h_2)$ or $h_3=\rev(h_4)$, $e$ and~$e'$ correspond to the same edge. In this case, we distinguish~$e$ and~$e'$ in the obvious way: when, for example, $h_1=\rev(h_2)$, the operation merely consists in suppressing the face~$f$ by gluing~$h_3$ and~$h_4$ together, the edge~$e$ and~$e'$ being distinguished in such a way that~$e$ is on the side where~$h_4$ was and~$e'$ is on the side where~$h_3$ was. See Figure~\ref{nn-}.

Remark that choosing the other possibility for the half-edges~$h_1$, $h_2$, $h_3$, $h_4$ merely changes the orientation of the loops considered and exchanges the left part with the right part. The final output, $(\q;e,e')$, remains the same. To convince oneself that this is true, one can think of the following similar situation: we cut a sphere in two halves, turn clockwise of some fixed angle one of the two hemispheres and glue the two parts together. The result is the same no matter what hemisphere we chose to rotate.

\bigskip

\begin{figure}[ht]
		\psfrag{x}[][][.8]{$h_1$}
		\psfrag{y}[][][.8]{$h_2$}
		\psfrag{v}[][][.8]{$v$}
		\psfrag{f}[][][.8]{$f$}
		\psfrag{l}[][][.8]{\textcolor{red}{$\pl$}}
		\psfrag{r}[][][.8]{\textcolor{red}{$\pr$}}
		\psfrag{e}[][][.8]{\textcolor{blue}{$e$}}
		\psfrag{g}[][][.8]{\textcolor{blue}{$e'$}}
		\psfrag{h}[][][.8]{\textcolor{red}{$e'$}}
		\psfrag{1}[][][.75]{$\pl_1$}
		\psfrag{2}[][][.75]{$\pl_2$}
		\psfrag{3}[][][.75]{$\pl_3$}
		\psfrag{4}[][][.75]{$\pl_4$}
		\psfrag{5}[][][.75]{$\pl_5$}
		\psfrag{6}[][][.75]{$\pl_6$}
		\psfrag{7}[][][.75]{$\pr_1$}
		\psfrag{8}[][][.75]{$\pr_2$}
		\psfrag{9}[][][.75]{$\pr_3$}
		\psfrag{a}[][][.75]{$\pr_4$}
		\psfrag{b}[][][.75]{$\pr_5$}
		\psfrag{c}[][][.75]{$\pr_6$}	
	\centering\includegraphics[width=.95\linewidth]{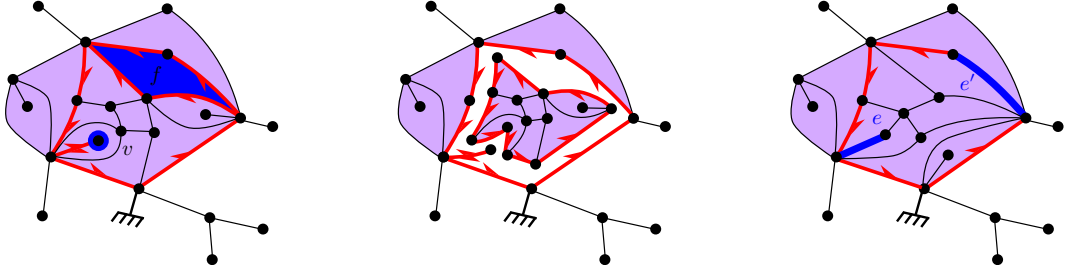}
	\caption{Removing a confluent face when the loops are pinched. Here, $d=3$.}
	\label{nn-2}
\end{figure}

Let us suppose now that~$\pp$ and~$\pp'$ meet before reaching~$v$. A path will be called a \emph{pinched loop} if it is of the form $\mathbbm a \ooo \mathbbm b \ooo \rev(\mathbbm b) \ooo \mathbbm c$, where~$\mathbbm b$ is a self-avoiding path and $\mathbbm a \ooo \mathbbm c$ is a simple loop that intersects~$\mathbbm b$ only at its origin. The paths $\rev(h_3) \ooo \pp \ooo \rev(\pp') \ooo \rev(h_4)$ and $h_1 \ooo h_2 \ooo \pp \ooo \rev(\pp')$ are now both pinched loops. We use the same conventions as above to define a left part and a right part, to glue them together, and to define~$e$ and~$e'$. The only difference is that we cut along the pinched part of the loop in the one of the two parts where it is possible, similarly as in Section~\ref{secpp}. See Figure~\ref{nn-2}.

In both cases, we set $m\de c$ ($c$ as confluent) and $\Phi_{n+1\dn,p}(\q';f,v) \de (\q;e,e',m)$.

\subsubsection{Simple face}\label{secnn-sf}

Let us now turn to the case where~$f$ is simple with respect to~$v$. The operation in this case is a little easier as it involves only one geodesic. We still let~$d$ be the distance between~$f$ and~$v$, and we let~$h_1$, $h_2$, $h_3$, $h_4$ be the four half-edges incident to~$f$ read in the counterclockwise order, $h_4^+$ being the closest to~$v$. We denote by~$\pp$ the left-most geodesic from~$h_4$ to~$v$.

\begin{figure}[ht]
		\psfrag{v}[][][.8]{$v$}
		\psfrag{f}[][][.8]{$f$}
		\psfrag{l}[][][.8]{\textcolor{red}{$\pl$}}
		\psfrag{r}[][][.8]{\textcolor{red}{$\pr$}}
		\psfrag{e}[][][.8]{\textcolor{blue}{$e$}}
		\psfrag{g}[][][.8]{\textcolor{blue}{$e'$}}
	\centering\includegraphics[width=.95\linewidth]{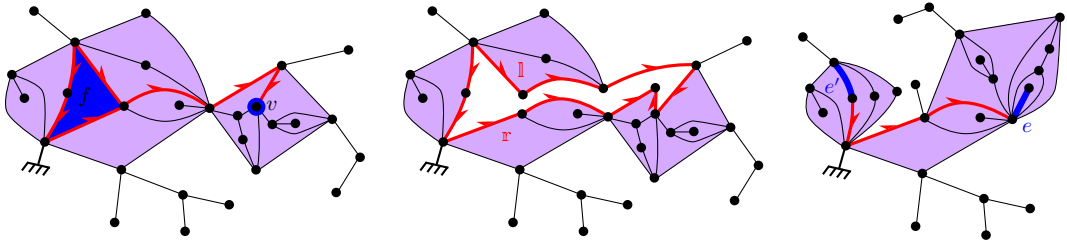}
	\caption{Removing a simple face. On this picture, $d=3$ and $m\de s_1$.}
	\label{nn-3}
\end{figure}

We cut along~$\pp$, starting from the corner delimited by~$h_4$ and $\rev(h_1)$ and stopping at~$v$. We obtain a map with a face of degree $2d+4$. In this map, every half-edge of~$\pp$ has two images, exactly one being incident to the face of degree $2d+4$. For $1\le k \le d$, let $\pp^r_k$ be the image of~$\pp_k$ that is incident to the face of degree $2d+4$, and~$\pp^l_k$ the other one. We set
$$\pl \de \rev(h_1)\ooo \pp^l \ooo \rev(\pp^r_d) \sand \pr \de h_2\ooo h_3 \ooo h_4 \ooo \pp^r_{1\to d-1}.$$
The map~$\q$ is then defined by gluing back~$\pl_k$ to~$\pr_k$ for $1\le k \le d+2$. The edge~$e$ is the one corresponding to~$\pl_{d+2}$ and~$e'$ the one corresponding to~$\pl_1$. See Figure~\ref{nn-3}.

\begin{figure}[ht]
		\psfrag{m}[][][.8]{$m\de s_1$}
		\psfrag{n}[][][.8]{$m\de s_2$}
		\psfrag{e}[][][.8]{$e$}
		\psfrag{f}[][][.8]{$e'$}
		\psfrag{1}[][][.8]{$h_1$}
		\psfrag{2}[][][.8]{$h_2$}
		\psfrag{3}[][][.8]{$h_3$}
		\psfrag{4}[][][.8]{$h_4$}
		\psfrag{v}[][][.8]{$v$}
		\psfrag{p}[][][.8]{$\pp$}
		\psfrag{q}[][][.8]{\textcolor{red}{$\pp'$}}
	\centering\includegraphics[width=.95\linewidth]{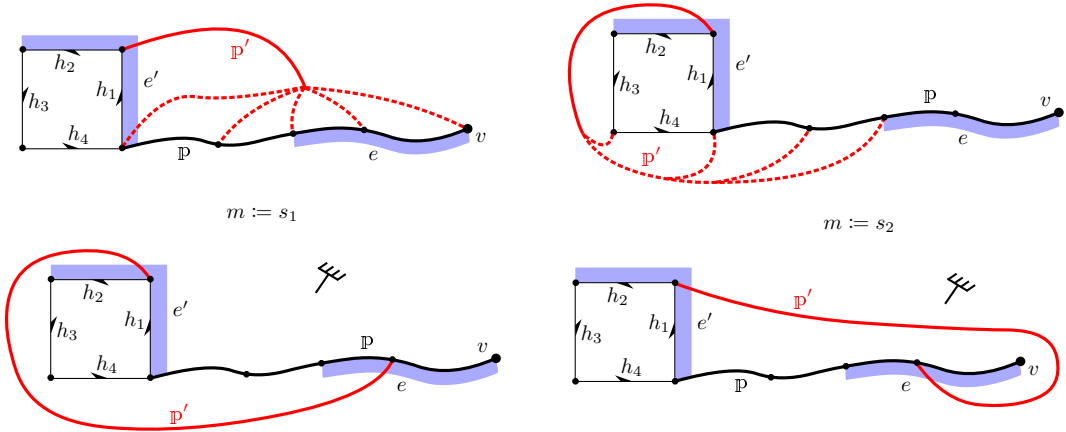}
	\caption{Setting the value of the mark~$m$. Beware that, on the first line, the external face is not necessarily in the infinite component of the plane whereas, on the second line, it is. The half-edges that will make up~$e$ and~$e'$ are highlighted.}
	\label{nn-4}
\end{figure}

The difficulty in this case is to set the value of the mark~$m$. To this end, we need to consider the left-most geodesic~$\pp'$ from~$h_1$ to~$v$ and see how it merges with~$\pp$ (see Figure~\ref{nn-4}). Let $k\in \{0,\dots,d+1\}$ be the smallest integer such that $\pp_k'^+ \in \{\pp_0^-,\dots,\pp_d^-,\pp_d^+\}$ (with the convention $\pp_0\de h_4$ and $\pp'_0\de h_1$). After cutting along~$\pp$,
\begin{itemize}
\item if $\pp_k'^+ \in \{\pp_1^{l-},\dots,\pp_d^{l-},\pp_d^{l+}\}$, we set $m\de s_1$;
\item if $\pp_k'^+ \in \{h_4^-,\pp_1^{r-},\dots,\pp_{d-1}^{r-}\}$, we set $m\de s_2$;
\item if $\pp_k'^+ = \pp_d^{r-}$ and the external face of~$\q'$ is to the right of the cycle $\pp'_{1\to d}\ooo \rev(\pp_{1\to d-1}) \ooo h_1$, then $m\de s_1$;
\item if $\pp_k'^+ = \pp_d^{r-}$ and the external face of~$\q'$ is to the left of $\pp'_{1\to d}\ooo \rev(\pp_{1\to d-1}) \ooo h_1$, then $m\de s_2$.
\end{itemize}
The reason why we choose this rule may seem surprising at this point but should become clear in the next sections. As a first insight, it can be noticed at this stage that~$h_1$ and~$h_2$ will make up~$e'$ and that $\pp_{d-1}^r$ and $\pp_{d}^r$ will make up~$e$. Here again, we set $\Phi_{n+1\dn,p}(\q';f,v) \de (\q;e,e',m)$.

\subsection{Adding a face}\label{secnn+}

We now present the mapping $\Phi_{n\up,p}$ from $\Qnp^{\textup{eem}}$ to $\Q_{n+1,p}^{\textup{fv}}$. Let $(\q;e,e',m) \in \Qnp^{\textup{eem}}$. Notice that the four distances between the extremities of~$e$ and the extremities of~$e'$ are either $d$, $d+1$, $d+1$, $d$ or $d$, $d+1$, $d+1$, $d+2$, for some $d\ge 0$. In the first case, we will say that~$e$ and~$e'$ are \emph{parallel}.

\begin{figure}[ht]
		\psfrag{d}[][]{$d$}
		\psfrag{+}[][]{$d+1$}
		\psfrag{2}[][]{$d+2$}
	\centering\includegraphics[width=.95\linewidth]{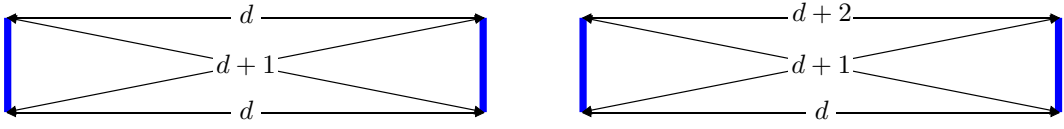}
	\caption{Parallel edges on the left, nonparallel edges on the right.}
	\label{}
\end{figure}

\subsubsection{Confluent case}

Let us start with the construction in the case $m=c$. In a first time, we moreover suppose that~$e$ and~$e'$ are parallel and we let~$d$ be the distance between them. We arbitrarily choose a half-edge~$\vec e$ corresponding to~$e$ and we consider the right-most geodesics~$\pp$ and~$\pp'$ from~$\vec e$ and from~$\rev(\vec e)$ toward~$e'$. As there are no cycles of odd length, it is easy to see that these paths do no intersect and, in particular, their endpoints are distinct: $\pp_d^+\neq \pp_d'^+$. Let~${\vec e}\,'$ denote the half-edge corresponding to~$e'$, directed from $\pp_d'^+$ to $\pp_d^+$. The loop
$${\vec e}\,' \ooo \rev(\pp) \ooo \rev(\vec e) \ooo \pp'$$
separates the map in a left and a right part as before, both of them having a face of degree $2d+2$. (If~$e$ and~$e'$ correspond to the same edge, one part will consist in a single edge and the other one will have a face of degree~$2$ where~$e$ and~$e'$ were). We denote by~$\pl$ and~$\pr$ the images of the loop in the left and in the right part. We define~$\q'$ by gluing the parts together while matching~$\pl_k$ with $\pr_{k+1}$ for $2\le k \le 2d +1$. This creates an extra face~$f$ incident to~$\pr_1$, $\pr_2$, $\rev(\pl_1)$, $\rev(\pl_{2d+2})$. We set $v\de \pr_{d+2}^+$ in the map~$\q'$ and $\Phi_{n\up,p}(\q;e,e',m)\de(\q';f,v)$. See Figure~\ref{nn-}.

Now, if~$e$ and~$e'$ are nonparallel, we set~$\vec e$ the half-edge corresponding to~$e$ directed toward~$e'$, as well as~$\vec e\,'$ the half-edge corresponding to~$e'$ directed away from~$e$. We also let $d-1$ be the distance between~$e$ and~$e'$. We consider the right-most geodesic~$\pp$ from~$\vec e$ toward $\vec e\,'^+$, and the right-most geodesic~$\pp'$ from~$\vec e$ toward~$\vec e\,'^-$. The path
$${\vec e}\,' \ooo \rev(\pp) \ooo \rev(\vec e) \ooo \vec e \ooo \pp'$$
is a pinched loop of length $2d+2$. (Note that, as soon as~$\pp$ and~$\pp'$ split, they cannot meet again.) As a result, it separates the map in a left and a right part as before and we use the same convention to define~$\q'$ and~$f$. The vertex~$v$ is defined as the endpoint of the pinched part in the part where the loop is pinched. We set $\Phi_{n\up,p}(\q;e,e',m) \de (\q';f,v)$.

\subsubsection{Simple case}\label{nn-sc}

\begin{figure}[ht]
		\psfrag{a}[][][.8]{$\vec e$}
		\psfrag{b}[][][.8]{$\vec e\,'$}
		\psfrag{e}[][][.8]{\textcolor{blue}{$e$}}
		\psfrag{f}[][][.8]{\textcolor{blue}{$e'$}}
		\psfrag{p}[][][.8]{$\pp$}
		\psfrag{m}[][]{$m=s_1$}
		\psfrag{n}[][]{$m=s_2$}		
	\centering\includegraphics[width=.95\linewidth]{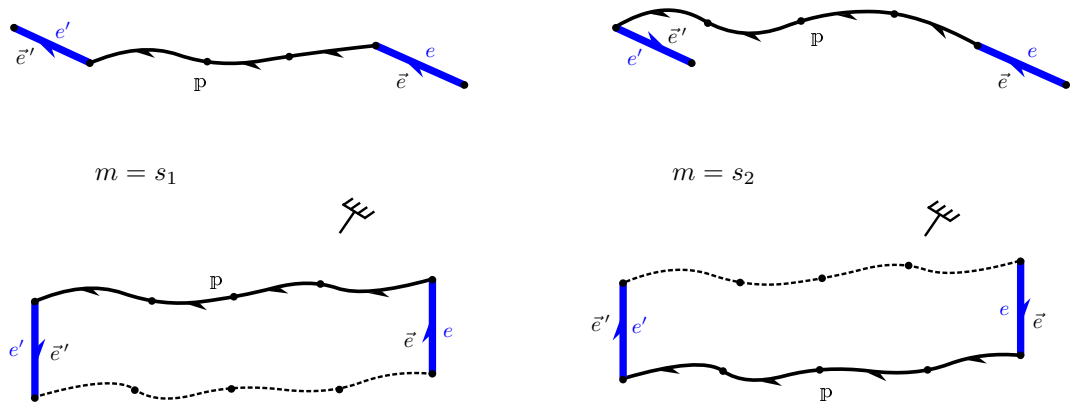}
	\caption{Choosing the path~$\pp$ when~$e$ and~$e'$ are nonparallel (top) or parallel (bottom) and when $m=s_1$ (left) or $m=s_2$ (right).}
	\label{pathp}
\end{figure}

We now suppose that $m\in\{s_1,s_2\}$. We will always proceed as follows: we will find a path~$\pp$ linking~$e$ to~$e'$, cut along it and create a new face by sliding the two sides of the cut as in Figure~\ref{nn-3}. Let us first describe how to choose the path~$\pp$ along which we will cut. See Figure~\ref{pathp}.

We first suppose that~$e$ and~$e'$ are nonparallel. Let~$\vec e$ be the half-edge corresponding to~$e$ directed toward~$e'$. If $m=s_1$, then~$\pp$ is the right-most geodesic from~$\vec e$ to~$e'$. If $m=s_2$, then~$\pp$ is the right-most geodesic from~$\vec e$ to the extremity of~$e'$ the farther away from~$e$.

Now, if~$e$ and~$e'$ are parallel, we consider the two right-most geodesics from both half-edges corresponding to~$e$ toward~$e'$. As there are no cycles of odd length, these two geodesics do not meet. As a result, it is possible to concatenate a half-edge corresponding to~$e$, one of these geodesics, a half-edge corresponding to~$e'$ and the reverse of the other geodesic so that the result is a simple loop. This may be done in two different ways, depending on which geodesic we choose to visit first, and choosing one way or the other merely changes the orientation of the loop. If $m=s_1$ (resp.\ $m=s_2$), we denote by~$\pp$ the geodesic such that if we visit it first, then the external face of the map lies to the right (resp.\ to the left) of this loop. 

In all cases, we denote by~$\vec e$ the half-edge corresponding to~$e$ directed toward the beginning of~$\pp$ and~$\vec e\,'$ the half-edge corresponding to~$e'$ directed away from the endpoint of~$\pp$. (We use an obvious extension in the case where~$\pp$ is the empty path.) We also set $d\de [\pp]$. By cutting along
$$\rev(\vec e \ooo \pp \ooo \vec e\,'),$$
we obtain a map with a face of degree $2d+4$. We denote by~$\pr$ the image of the previous path that is incident to this face, and~$\pl$ the image that is not incident to this face. Beware that in the case where~$e$ and~$e'$ are nonparallel and $m=s_2$, one of these paths is not self-avoiding if $\pp_d^-=\vec e\,'^+$. In this case, the part of the map separated by the loop $\vec e\,' \ooo \pp_d$ that does not contain~$e$ remains attached to~$\vec e\,'^+$. We then define~$\q'$ by gluing back~$\pl_k$ with $\pr_{k+2}$ for $2\le k \le d$, as well as $\pl_{d+1}$ with $\rev(\pl_{d+2})$. This creates an extra face~$f$ incident to~$\pr_1$, $\pr_2$, $\pr_3$ and $\rev(\pl_1)$. We set $v\de \pl_{d+1}^+$ in the map~$\q'$ and $(\q';f,v) \de \Phi_{n\up,p}(\q;e,e',m)$. See Figure~\ref{nn-3}.

\subsection{These mappings are reverse one from another}

\begin{thm}\label{thmnn}
The mappings $\Phi_{n\up,p}:\Qnp^{\textup{eem}} \to \Q_{n+1,p}^{\textup{fv}}$ and $\Phi_{n+1\dn,p}:\Q_{n+1,p}^{\textup{fv}} \to \Qnp^{\textup{eem}}$ are one-to-one and reverse one from another.
\end{thm}

\begin{pre}
We will proceed in a similar manner as in the proof of Theorem~\ref{thmpp}. Our constructions from Sections~\ref{secnn-} and~\ref{secnn+} are clearly reverse one from another, the only thing we have to check is that the paths along which we cut correspond. In other words, we need to verify that the path denoted by~$\pl$ in both sections is the same if we apply our constructions to a map and to its image through one of our mappings, and that the same goes for the path~$\pr$.

\paragraph{Confluent case.}

Let $(\q';f,v) \in \Q_{n+1,p}^{\textup{fv}}$ be such that~$f$ is confluent with respect to~$v$ and let $(\q;e,e',m) \de \Phi_{n+1\dn,p}(\q';f,v)$. We use here the notation of Section~\ref{secnn-}. Let us first suppose that the loop considered in Section~\ref{secnn-cf} is not pinched (case of Figure~\ref{nn-}). In order to conclude that $\Phi_{n\up,p}(\q;e,e',m) = (\q';f,v)$, it will be sufficient to show that $\rev(\pr_{2\to d+1})$ and $\pr_{d+3\to 2d+2}$ are the two right-most geodesics from the two half-edges corresponding to~$e$ toward~$e'$. Note that, in particular, this will entail that~$e$ and~$e'$ are parallel. This fact is again shown by the winding argument: if, for example, the right-most geodesic from $\rev(\pr_{d+2})$ to~$e'$ leaves $\rev(\pr_{2\to d+1})$ and enters the left part or the right part, then it has to leave it or it would contradict either the fact that~$\pp$ and~$\pp'$ are left-most geodesics from~$f$ to~$v$ or the fact that~$h_1^+$ and~$h_3^+$ are at distance $d+1$ from~$v$. Repeating the argument, it will indefinitely wind around~$e'$.

If the loop of Section~\ref{secnn-cf} is pinched (case of Figure~\ref{nn-2}), we conclude by exactly the same argument. Note that in this case, $e$ and~$e'$ are nonparallel.

\paragraph{Simple case.}

Let us suppose now that $(\q';f,v) \in \Q_{n+1,p}^{\textup{fv}}$ is such that~$f$ is simple with respect to~$v$ and let again $(\q;e,e',m) \de \Phi_{n+1\dn,p}(\q';f,v)$. We use the notation of Section~\ref{secnn-sf} (see Figure~\ref{vecevece}). In~$\q$, let us set $\vec e \de \rev(\pl_{d+2})$ and $\vec e\,' \de \rev(\pl_1)$. It is not hard (although it requires some care) to show that $\rev(\pl_{2\to d+1})$ is the right-most geodesic from~$\vec e$ to~$\vec e\,'^-$. This entails in particular that $d_\q\lp \vec e\,^+,\vec e\,'^- \rp = d$. It is also easy to see that $d_\q\lp \vec e\,^-,\vec e\,'^- \rp = d+1$.

\begin{figure}[ht]
		\psfrag{a}[][][.8]{$\vec e$}
		\psfrag{b}[][][.8]{$\rev(\pl_{2\to d+1})$}
		\psfrag{c}[][][.8]{$\vec e\,'$}
		\psfrag{e}[][][.8]{$e$}
		\psfrag{f}[][][.8]{$e'$}
		\psfrag{1}[][][.8]{$h_1$}
		\psfrag{2}[][][.8]{$h_2$}
		\psfrag{3}[][][.8]{$h_3$}
		\psfrag{4}[][][.8]{$h_4$}
		\psfrag{v}[][][.8]{$v$}
		\psfrag{p}[][][.8]{}	
		\psfrag{r}[][][.8]{$\pr$}	
		\psfrag{l}[][][.8]{\textcolor{red}{$\pl$}}		
	\centering\includegraphics[width=.95\linewidth]{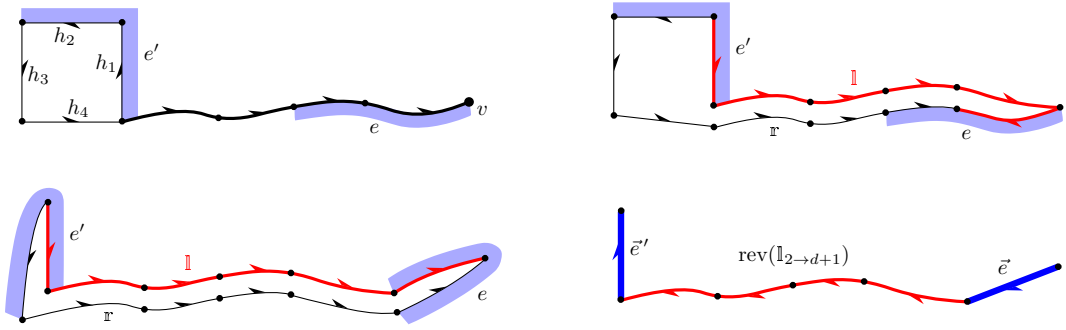}
	\caption[Reminder of the construction in the case of a simple face and some notation.]{Reminder of the construction $\Phi_{n+1\dn,p}$ in the case of a simple face and some notation.}
	\label{vecevece}
\end{figure}

The remaining pair of distances
$$\lp d_\q\lp \vec e\,^+,\vec e\,'^+ \rp, d_\q\lp \vec e\,^-,\vec e\,'^+ \rp \rp$$
may be either $(d+1,d+2)$, $(d-1,d)$ or $(d+1,d)$, depending on how~$\pp'$ merges with~$\pp$. The first two cases correspond to the top line of Figure~\ref{nn-4}. In these cases, $e$ and~$e'$ are nonparallel and~$\vec e$ is directed toward~$e'$. As a result, the path $\rev(\pl_{2\to d+1})$ corresponds to the path of Section~\ref{nn-sc}. The last case corresponds to the bottom line of Figure~\ref{nn-4}; we obtain that~$e$ and~$e'$ are parallel, and distinguishing whether $m=s_1$ or $m=s_2$, we see that $\rev(\pl_{2\to d+1})$ still corresponds to the path of Section~\ref{nn-sc}. This is suficient to conclude that $\Phi_{n\up,p}(\q;e,e',m) = (\q';f,v)$.

\bigskip

The fact that $\Phi_{n+1\dn,p} \circ \Phi_{n\up,p}$ is the identity over $\Qnp^{\textup{eem}}$ also follows from the same kind of arguments and is left to the reader.
\end{pre}

\subsection{Interpretation through the \BDG bijection}\label{secbdg}

The \BDG bijection is a classical bijection allowing to encode maps with simper objects. In our case of quadrangulations, it specializes into a bijection between quadrangulations having a distinguished vertex and so-called well-labeled forests. We will not describe this bijection in this work and refer to~\cite{bouttier04pml} for a complete description. See also~\cite{bettinelli11slr,bouttier09dsqb} for an exposition in the particular case of quadrangulations.

\begin{defi}
A \emph{well-labeled forest} is a pair $(\f,\lab)$ where~$\f=(\tr_1,\dots,\tr_p)$ is a forest and $\lab$ is an integer-valued function on the vertices of~$\f$ satisfying the following: 
\begin{itemize}
	\item $|\lab(u) - \lab(v)| \le 1$ whenever~$u$ and~$v$ are vertices of the same tree linked by an edge,
	\item $\lab(\rho_{i+1}) \ge \lab(\rho_{i})-1$ for all $1\le i \le p$, where~$\rho_i$ denotes the root vertex of~$\tr_i$, and $\rho_{p+1} \de \rho_1$,
	\item $\lab(\rho_{1})=0$.
\end{itemize}
\end{defi}

Our quadrangulations in $\Q_{n+1,p}^{\textup{fv}}$ come with a distinguished vertex, so it seems natural to try to understand what happens to the coding well-labeled forest through $\Phi_{n+1\dn,p}$. However, after applying our bijection, the forest is no longer a forest, it becomes some map with two faces with rather complicated rules on its labels. We may also mention that the maps in $\Qnp^{\textup{eem}}$ do not come a priori with a natural distinguished vertex.

\bigskip

A quite remarkable fact that follows from our theorems is that a third of the quadrangulations in $\Q_{n+1,p}^{\textup{fv}}$ are such that the distinguished face is confluent with respect to the distinguished vertex, one third are organized as on the left of Figure~\ref{nn-4} and one third are organized as on the right of Figure~\ref{nn-4}. This fact is actually very easy to show by using the \BDG bijection.

In fact, quadrangulations in $\Q_{n+1,p}^{\textup{fv}}$ correspond to well-labeled forests having a distinguished edge. The first third correspond to the case where the edge links two vertices having the same label. In the other cases, the edge links a vertex~$v$ labeled~$\ell$ to a vertex~$v'$ labeled $\ell+1$. Removing the edge breaks one of the trees~$\tr$ in two parts. We consider, on the one hand, the part of~$\tr$ not containing its root and, on the other hand, the part of~$\tr$ containing the root, together with all the other trees. We denote by~$S_v$ the one of these two sets containing~$v$ and $S_{v'}$ the one containing~$v'$. Finally, we set $m_v \de \min_{S_v} \lab - \ell$ and $m_{v'} \de \min_{S_{v'}} \lab - (\ell+1)$. Then, the top-left part of Figure~\ref{nn-4} corresponds to the case where $m_{v'} < m_{v}$, the top-right part corresponds to the case where $m_{v} < m_{v'}$. The bottom line corresponds to the case where $m_{v} = m_{v'}$; on the left, the root of~$\tr$ belongs to~$S_v$, on the right, it belongs to~$S_{v'}$. By symmetry, we recover the distribution into three thirds.


\bibliographystyle{alpha}
\bibliography{__these_bib}
\end{document}